\documentstyle[12pt,fleqn,leqno]{article}
\newcommand{\be}{\begin{equation}}
\newcommand{\ee}{\end{equation}}
\newcommand{\ba}{\begin{array}}
\newcommand{\ea}{\end{array}}
\newcommand{\bi}[1]{
\renewcommand{\arraystretch}{0.4}
(\!\! \ba{c} \scriptscriptstyle {#1} \\ \scriptscriptstyle  2\ea\!\!)
\renewcommand{\arraystretch}{1}
}
\newcommand{\r}{Ramanujan}
\renewcommand{\a}{\alpha}
\renewcommand{\l}{\lambda}
\newcommand{\re}{\mbox{Re$\;$}}
\newcommand{\im}{\mbox{Im$\;$}}
\newcommand{\mg}{\rm}
\renewcommand{\em}{\it}
\newcommand{\bea}{\begin{eqnarray}}
\newcommand{\eea}{\end{eqnarray}}

\newcommand{\dis}{\vspace{-0.5\abovedisplayskip}}
\newcommand{\Int}{\int_{-\infty}^\infty}
\newcommand{\Sum}{\sum_{n=-\infty}^\infty}
\newcommand{\Summ}{\sum_{n=0}^\infty}
\newcommand{\Inte}{\int_{0}^\infty}
\newcommand{\Prod}{\prod_{n=0}^\infty}
\newcommand{\Integ}{\int_{0}^1}

\newcommand{\ph}[2]{\mbox{${}_{#1}\psi_{#2} $}}
\newcommand{\phy}[2]{\mbox{${}_{#1}\phi_{#2} $}}
\newcommand{\g}[1]{\Gamma({#1})}
\newcommand{\G}[1]{{\Gamma_q({#1})}}
\newcommand{\Ga}[1]{{\Gamma_{q^2}({#1})}}
\begin{document}
\title{\bf Some Basic Bilateral \\Sums and Integrals\thanks{
This work was supported, in part, by NSF grant DMS 9203659 and NSERC grant \# A6197.} }
\author{Mourad E. H.  Ismail\thanks{Department  of Mathematics, University of
South Florida, Tampa, FL 33620, U. S. A.} and Mizan Rahman\thanks{
Department of Mathematics and statistics, Carleton University, Ottawa, Ontario,
Canada, K1S 5B6.} }
\date{}
\maketitle
\begin{abstract}
By splitting the real line into intervals of unit length a doubly infinite
integral of the form $\Int F(q^x)\,dx,\; 0<q<1$, can clearly be expressed
as $\Integ \Sum F(q^{x+n})\,dx$, provided $F$ satisfies the appropriate
conditions. This simple idea is used to prove \r 's integral analogues of
his \ph{1}{1} sum and give a new proof of Askey and Roy's extention
of it. Integral analogues of the well-poised \ph{2}{2} sum as well as the
very-well-poised \ph{6}{6} sum are also found in a straightforward manner.
An extension to a very-well-poised and balanced \ph{8}{8} series is also
given. A direct proof of a recent q-beta integral of Ismail and Masson 
is given.
\end{abstract}

\bigskip
{\bf Running title}:$\;$ Sums and Integrals

\bigskip
{\em 1990  Mathematics Subject Classification}: Primary 33D05, Secondary 33D20

{\em Key words and phrases}. Q-beta integrals, Askey-Roy and \r\ integrals, evaluation of sums and integrals, weight functions, $_6\psi_6$ sum.
\vfill\eject
\section{ Introduction.} The familiar form of the classical beta integral of
Euler is
\be
 B(a,b)=\Integ t^{a-1}(1-t)^{b-1}\,dt=\frac{\g{a}\,\g{b}}{\g{a+b}}\,,\label{b1}
\ee
\re $(a, b)>0$. A less familiar form, obtained by a simple change of variable,
is
\be 
B(a,b)=\Inte \frac{t^{a-1}\,dt}{(1+t)^{a+b}}\,.\label{b2} 
\ee
There have been many extensions of both these forms, see, for example,
Askey [2--5], Askey and Roy [6], Gasper [9, 10], Rahman and Suslov [18], and
the references therein. A ``curious'' extension of (\ref{b2}) that was given
by \r\ [21] in 1915 is
\be
 \Inte t^{a-1}\,\frac{(-tq^{a+b};\,q)_\infty}{(-t;\,q)_\infty}\,dt=\frac
{\g{a}\,\g{1-a}}{\G{a}\,\G{1-a}}\,\frac{\G{a}\,\G{b}}{\G{a+b}}\,,\label{b3}
\ee
where \re $(a,b)>0,\; 0<q<1$, the $q$-gamma function $\G{x}$ is defined by
\be
\G{x}=\frac{(q;\,q)_\infty}{(q^x;\,q)_\infty}\,(1-q)^{1-x}\,,\quad x\neq 0,-1,
-2,\ldots,\label{b4} 
\ee
and the infinite products by
\be
(a;\,q)_\infty=\Prod(1-aq^n)\,.\label{b5} 
\ee
The fact that the limit of the formula (\ref{b3}) as $q\to1^-$ is (\ref{b2}) follows from the
properties
\be
\lim_{q\to1^-}\G{x}=\Gamma(x),\quad \lim_{q\to1^-}\frac{(-tq^c;\,q)_\infty}{
(-t;\:q)_\infty}=\frac1{(1+t)^c},\label{b6}
\ee
see [11].\par

Askey and Roy [6] introduced a third parameter into the formulas and gave
the following extension of (\ref{b3})
\be
\Inte t^{c-1}\,\frac{(-tq^{b+c},-q^{a-c+1}/t;\,q)_\infty}{
(-t,-q/t;\,q)_\infty}\,dt=
\frac{\g{c}\,\g{1-c}}{\G{c}\,\G{1-c}}\,\frac{\G{a}\,\G{b}}{\G{a+b}}\,,\label{b7} 
\ee
which holds for \re$(a,b,c)>0$. In the limit $c\to0^+$ this becomes
\be
\Inte \frac{(-tq^{b},-q^{a+1}/t;\,q)_\infty}{
(-t,-q/t;\,q)_\infty}\,\frac{dt}{t}=
\frac{\log{q^{-1}}}{1-q}\,\frac{\G{a}\,\G{b}}{\G{a+b}}\,,\label{b8} 
\ee
which restores the symmetry in $a$ and $b$ that was there in
both (\ref{b1}) and (\ref{b2}), but not in (\ref{b3}), see Gasper [9--10].
Following [11] we have used the  contracted notation
\be
(a_1,a_2,\ldots,a_k;\,q)_\infty=\prod_{j=1}^k (a_j;\,q)_\infty\,.\label{b9} 
\ee
\par

Hardy [12] gave a proof of (\ref{b3}) that \r\ did not, and discussed \r's
general method of evaluating such integrals in [13]. Askey [2] gave another
proof of (\ref{b3}). Askey's method is rather close to the Pearson-type 
first order difference equation technique that has been used extensively
by the Russian school of Nikiforov, Suslov and Uvarov, see for example,
[17, 23], as well as of Atakishiyev and Suslov [7]. It was pointed out in
[18] and [20] that the origin of both Barnes and \r-type integrals can be
traced to a Pearson equation on linear, $q$-linear, quadratic or
 $q$-quadratic
lattice with appropriately chosen coefficient functions so that the boundary
conditions can be satisfied in the two cases. In [20] Rahman and Suslov
found what they consider a better way of dealing with the \r-type integrals,
and evaluated extensions of some of \r-type formulas as well as extensions
of the summation formulas of Gauss, and, Pfaff and Saalsch\"{u}tz. 
The idea is very simple.
Suppose that $f(x)$ is continuous on $[a,\infty)$, has no singularities on 
the real line, and its integral on $[a,\infty)$ exists. Suppose also that
$\Sum f(x+n)$ converges uniformly for $x\in [a,a+1]$. Then
\be
\int_a^\infty f(x)\,dx=\int_a^{a+1}\Summ f(x+n)\,dx\,.\label{b10}
\ee
For integrals on the whole real line the corresponding formula is
\be
\Int f(x)\,dx=\Integ\Sum f(x+n)\,dx\,,\label{b11}
\ee
provided, of course, that the bilateral sum $\sum_{-\infty}^\infty f(x+n)$
converges uniformly for $x\in[0,1]$. What this method does is to establish a
direct correspondence between the integrals on the left side of 
(\ref{b10})
and (\ref{b11}), and the infinite series on the right. 
So  for the method to be useful
we have to be able to handle the infinite sums so that we can apply this 
knowledge to compute the infinite integrals. As is well-known in classical 
analysis, it is often the case that an infinite integral over a function is
easier to compute than an infinite sum. So the method described above has
a very limited applicability. It is applicable when the series inside the
integrals on the right sides of (\ref{b10}) and (\ref{b11}) are summable
(meaning that the sum can be evaluated in closed forms) or at least 
transformable in a way that the ensuing formulas are simpler. Such is the
case for some hypergeometric and basic hypergeometric series, bilateral or
otherwise. \par

A basic bilateral series in base $q$ (assumed throughout this paper to 
satisfy $0<q<1$), with $r$ numerator and $r$ denominator parameters is
defined by
\be
\ph{r}{r}\left[\ba{c}
a_1,a_2,\ldots,a_r\\
b_1,b_2,\ldots,b_r \ea ;\,q,z\right]=\, _r\psi_r(a_1, \cdots, a_r; b_1, \cdots, b_r; q, z) =\sum_{n=-\infty}^\infty
\frac{(a_1,a_2,\ldots,a_r;\,q)_n}{(b_1,b_2,\ldots,b_r ;\,q)_n}\,z^n\,,\label{b12}
\ee
where\dis
\[
(a_1,a_2,\ldots,a_r;\,q)_n=\prod_{j=1}^r(a_j;\,q)_n\,,
\]\dis
\be
(a_j;\,q)_n=\left\{\ba{ll}
1,  & \mbox{if }n=0,\\
\prod_{k=0}^{n-1}(1-a_jq^k), & \mbox{if }n=1,2,\ldots. \ea \right.\label{b13}
\ee
The series (\ref{b12}) is absolutely convergent in the annulus
\be
\left| \frac{b_1b_2\ldots b_r}{a_1a_2\ldots a_r}\right| < |z|<1.\label{b14}
\ee
If any  one of the denominator parameters equals $q$, say, $b_r=q$,
then the first non-zero term in the series corresponds to $n=0$, and
the series becomes a basic generalized hypergeometric series:
\bea\lefteqn{
\phy{r}{r-1}\left[\ba{c}
a_1,a_2,\ldots,a_r\\
b_1,b_2,\ldots,b_{r-1} \ea ;\,q,z\right]= \mbox{}_r\phi_{r-1}(a_1,\cdots, a_r; b_1, \cdots, b_{r-1}; q, z)\label{b15}}\\
& & := \Summ
\frac{(a_1,a_2,\ldots,a_r;\,q)_n}{(q, b_1,b_2,\ldots,b_{r-1} ;\,q)_n}\,z^n, \nonumber
\eea
which is absolutely convergent inside the unit circle $|z|=1$, for
further details see [11]. \par

One of the most important evaluations of a basic bilateral hypergeometric
series is the one due to \r\ [13]
\be
\ph{1}{1}\left[\ba{c} a\\b \ea;\,q,z\right]=
\sum_{n=-\infty}^\infty\frac{(a;\,q)_n}{b;\,q)_n}\,z^n 
=\frac{(q,b/a,az, q/az;\,q)_\infty}{(b,q/a,z,b/az;\,q)_\infty}\,.\label{b16}
\ee
Many different proofs of this formula have appeared in the literature, but
the ones that are most often quoted and instructive are in [1] and [14].
 However,
one runs into trouble with a bilateral series immediately
 after the \ph{1}{1} level. Instead of a nice compact
 formula like the $q$-Gauss formula
\be
\phy{2}{1}\left[\ba{r}a,\;b\\c\ea;\,q,c/ab\right]=\frac{(c/a,c/b;\,q)_\infty}{
(c,c/ab;\,q)_\infty}\,,\quad|c/ab|<1, \label{b17}
\ee
one has a 2-term formula for the corresponding \ph{2}{2} sum:
\bea
\lefteqn{\ph{2}{2}\left[\ba{ll}a,&b\\c,&d\ea;\,q, \; cd/abq\right]}\label{b18}\\
 & & \mbox{}-\frac{\a}{q}\,
\frac{(q/c,q/d,\a/a,\a/b;\,q)_\infty}{(q/a,q/b,\a/c,\a/d;\,q)_\infty}\:
\ph{2}{2}\left[\ba{ll}aq/\a,&bq/\a\\cq/\a,&dq/\a\ea;\,q,cd/abq\right] \nonumber\\
 & & = \frac{(\a,q/\a,cd/\a q,\a q^2/cd,q,c/a,c/b,d/a,d/b;\,q)_\infty}{(c/\a,\a q/c,d/\a,\a q/d, c,d,q/a,q/b,cd/abq;\,
q)_\infty}, \nonumber
\eea
see [11], where $\a$ is an arbitary parameter such that no zeros appear in
the denominators. It is clear that this formula reduces to (\ref{b17}) when
$d=q$. This is not a very well-known formula but a special case of it was
mentioned in [3]. As the number of parameters of the summation
formulas increases, one needs to impose more restrictions.
 The bilateral series that has the most desirable structure
is the very-well-poised one, namely
\bea
\lefteqn{\ph{r+2}{r+2}\left[ \ba{c}qa^{1/2},-qa^{1/2},\;\;a_1,\;\;a_2,\;\;\ldots,\;\;\;a_r \\
a^{1/2},-a^{1/2},qa/a_1,qa/a_2,\ldots,qa/a_r\ea\:;q\:,z\right]}\label{b19}\\
& & =\sum_{n=-\infty}^\infty\frac{(1-aq^{2n})(a_1,a_2,\ldots,a_r ;\,q)_n}
{(1-a)(qa/a_1,qa/a_2,\ldots,qa/a_r;q)_n}\:z^n\,.\nonumber 
\eea
The most general summation formula for a basic bilateral series is
Bailey's [8] \ph{6}{6} sum:
\bea
\lefteqn{\ph{6}{6}\left[ 
\ba{c}qa^{1/2},-qa^{1/2}\;,\;\;b,\;\;\;c,\;\;\;d,\;\;\;e \\
a^{1/2},-a^{1/2},aq/b,aq/c,aq/d,aq/e \ea\,;\:q,qa^2/bcde\right]}\label{b20} \\
& & \mbox{}=
\frac{(q,q/a,aq,aq/bc,aq/bd,aq/be,aq/cd,aq/ce,aq/de;\:q)_\infty}{
(q/b,q/c,q/d,q/e,aq/b,aq/c,aq/d,aq/e,qa^2/bcde\,;\:q)_\infty}\,,\nonumber
\eea
provided that $|qa^2/bcde|<1$, see [11]. Clearly, a function $f(x)$ for which
$\sum_{n=-\infty}^\infty f(x+n)$ corresponds to the sum on the left side of
(\ref{b20}) has to be of interest as far as the applicability of (\ref{b11})
is concerned. Accordingly, we first rewrite this formula in the form:
\bea\lefteqn{\sum_{n=-\infty}^{\infty}(aq^{n+1}/b,aq^{n+1}/c,aq^{n+1}/d,
aq^{n+1}/e;\:q)_\infty}\label{b21} \\
& & \cdot(q^{1-n}/b,q^{1-n}/c,q^{1-n}/d,q^{1-n}/e;\:q)_\infty 
\cdot(1-aq^{2n}) a^{2n}q^{n^2-n}\nonumber\\
& & =\frac{(q,a,q/a,aq/bc,aq/bd,aq/be,aq/cd,aq/ce,aq/de;\:q)_\infty}{
(qa^2/bcde;\; q)_\infty}\,.\nonumber
\eea
This suggests considering an integral of the form
\bea 
\lefteqn{ J:=\Int(aq^{x+1}\!/b,aq^{x+1}\!/c,aq^{x+1}\!/d, aq^{x+1}\!/e,
q^{1-x}\!/b,q^{1-x}\!/c,q^{1-x}\!/d;\:q)_\infty}\label{b22} \\ 
& & \times (q^{1-x}\!/e;\:q)_\infty (1-aq^{2x}) a^{2x}q^{x^2-x} \omega(x)
\,dx \nonumber
\eea
where $\omega(x)$ is a bounded continuous unit-period function on {\bf R},
i.e, $\omega(x\pm1)=\omega(x)$. \par

We shall evaluate this integral, (1.22), in \S3 by using (\ref{b11}) and (\ref{b21}),
and consider an extension of it in \S5. In \S2, however, we shall deal with
an integral analogue of (\ref{b16}) essentially showing that \r's formula
(\ref{b3}) is precisely that analogue. For an integral analogue of (\ref{b18})
we refer the reader to [20]. As a straightforward application of (\ref{b8}) 
we will also show in \S2 how to obtain a $q$-analogue of an integral of \r\
involving a product of two Bessel functions where the variable is the order
and not the argument. In \S4 we shall consider a case where $\omega(x)$ is a
unit anti-periodic function in (\ref{b22}), i.e $\omega(x\pm1)=-\omega(x)$,
thereby obtaining a $q$-analogue of yet another formula due to \r.

Ismail and Masson [16] proved that if the $q^{-1}$-Hermite polynomials are orthogonal with respect to a probability measure $d\psi$ then 
\bea\lefteqn{
\int^\infty_{-\infty}\prod^4_{j=1} (-t_j (x+\sqrt{x^2+1}),t_j(\sqrt{x^2+1}-x);
 q)_\infty\, d\psi (x)} \label{c23} \\
& & =\left[\prod_{1\leq j < k\leq 4} (-
t_jt_k/q; q)_\infty\right]/(t_1t_2t_3t_4q^{-3}; q)_\infty,\nonumber
\eea
holds, provided that the integral exists.  The corresponding moment 
problem has infinitely many solutions so one expects (1.23) to lead 
to an evaluation of an infinite family of integrals. Ismail and Masson 
[16] pointed out that Bailey's $_6\psi_6$ sum is (1.23)  with $d\psi$ a
 general extremal measure of the $q^{-1}$-Hermite moment problem. 
Ismail and Masson also observed that Askey's $q$-beta integral (3.4) 
corresponds to an absolutely continuous $d\psi$ and in this sense
 (3.4) is a continuous analogue of the $_6\psi_6$ sum of (1.20). Ismail 
and Masson [16] proved that the $q^{-1}$-Hermite polynomials are 
orthogonal with respect to the absolutely continuous measure $d\mu(x;\eta)$,
 where
\begin{equation}
\frac{d\mu(x;\eta)}{dx} = \frac{e^{2\eta_1}\sin\eta_2 \cosh \eta_1(q,\,
-qe^{2\eta_1}, -qe^{-2\eta_1};q)_\infty |(qe^{2i\eta_2};q)_\infty |^2}
{\pi\; |(e^{\xi+\eta},
\,-e^{\eta - \xi},\, -qe^{\xi-\eta},\, qe^{-\xi-\eta};\,q)_\infty|^2}, 
\; x = \sinh \xi, \eta = \eta_1 + i\eta_2.
\end{equation}
This and (1.23) led Ismail and Masson to the q-beta integral
\begin{equation}
\int^\infty_{-\infty}\frac{\prod^4_{j=1} (-t_je^{\xi},t_je^{-\xi};q)_\infty } {|(e^{\xi+\eta},
\,-e^{\eta - \xi},\, -qe^{\xi-\eta},\, qe^{-\xi-\eta};q)_\infty|^2}\,\cosh \xi\; d\xi \quad \quad
\end{equation}
$$\quad =\frac{\pi e^{-2\eta_1}\prod_{1\leq j < k\leq 4} (-
t_jt_k/q;q)_\infty}{\sin\eta_2 \cosh \eta_1(q, t_1t_2t_3t_4q^{-3},\,
-qe^{2\eta_1}, -qe^{-2\eta_1};q)_\infty |(qe^{2i\eta_2};q)_\infty |^2}.
$$
In Section 3 we shall give a direct proof of (1.25) and show that this
q-beta integral is another continuous analogue og the $_6\psi_6$ sum.

\section{ The Askey-Roy Integral and an Application.}
\setcounter{equation}{0} 
In this section we give new evaluations of (1.7) and (1.8) and also give some applications of them. Let us first rewrite (\ref{b16}) in the form
\renewcommand{\arraystretch}{0.4}
\bea
\lefteqn{
\sum_{n=-\infty}^\infty(bq^n,q^{1-n}\!/a;\:q)_\infty\,q^{(\!\!\ba{c}
\scriptscriptstyle n\\ \scriptscriptstyle 2\ea\!\!)
}(az)^ne^{\pi i n}}\label{c1} \\
 & & =\frac{(q,b/a,az,q/az;\:q)_\infty}{(z,b/az;\:q)_\infty}\,
,\;\;|b/a|<|z|<1.\nonumber
\eea
\renewcommand{\arraystretch}{1}
This suggests that we consider the integral
\be
I:=\Int (bq^x,q^{1-x}\!/a;\:q)_\infty\,q^{\bi{x}}(az)^x\,e^{\pi i x}\,
\omega(x)\, dx\,,\label{c2}
\ee
where $\omega(x\pm1)=\omega(x)$. For $0<q<1,\;|b/a|<|z|<1$ and continuous
bounded functions $\omega(x)$, it is clear that the integral exists.

{\bf Proof of (1.8)}. By (\ref{b11}) and (\ref{b16}) we have
\bea
\lefteqn{I=\Integ(bq^x,q^{1-x}\!/a;\:q)_\infty\,q^{\bi{x}}(az)^x\,e^{\pi i x}\,
\omega(x)\,\ph{1}{1}\left[\ba{r}aq^x\\bq^x\ea;\;q,z\right]\,dx} \label{c3}\\
& & =\frac{(q,b/a;\:q)_\infty}{(z,b/az;\:q)_\infty}\,\Integ(azq^x,
q^{1-x}\!/az;\:q)_\infty\,q^{\bi{x}}(az)^x\,e^{\pi i x}\,\omega(x)\, dx\,.
\nonumber
\eea
It is easy to see that the integrand in the last line above is unit periodic,
so it can be absorbed in $\omega$. Set
\be
\omega(x)=\frac{p(x)}
{(azq^x, q^{1-x}\!/az;\:q)_\infty\,q^{\bi{x}}(az)^x\,e^{\pi i x}}\,,\quad 
p(x\pm1)=p(x)\,, \label{c4}
\ee
which when substituted in (\ref{c2}) and (\ref{c3}), gives the formula
\be
\Int\frac{(bq^x,q^{1-x}\!/a;\:q)_\infty}{(azq^x,q^{1-x}\!/az;\:q)_\infty} p(x)\, dx
=\frac{(q,b/a;\:q)_\infty}{(z,b/az;\:q)_\infty}\,\Integ p(x)\,dx.\label{c5}
\ee
Assuming that $p(x)$ is independent of $a,b$ and $z$, let us now
replace $a,b,z$ by $-q^{-a},-q^b,-q^a$ in (\ref{c5}) to get
\bea
\lefteqn{ \Int  \frac{(-q^{b+x},-q^{a+1-x};
q)_\infty}{(-q^x,-q^{1-x};q)_\infty}p(x)dx }\label{c6} \\
& & =\frac{(q,q^{a+b};\:q)_\infty}{(q^a,q^b;\:q)_\infty}\,\Integ p(x)\,dx
=(1-q)^{-1}\frac{\G{a}\,\G{b}}{\G{a+b}}\,\Integ p(x)\,dx. \nonumber
\eea
Setting $p(x)\equiv 1$ and changing the variable by $q^x=t$, we establish 
(\ref{b8}) and the proof is complete. \par

The proof of (1.7)  is based on a different choice of $\omega$. Instead of (\ref{c4}), choose
\be
\omega(x)=\frac{q^{-\bi{x}}\,p(x)}{(-q^x,-q^{1-x};\:q)_\infty}\,,\quad 
p(x\pm1)=p(x)\,. \label{c7}
\ee

{\bf Proof of (1.7)}.  The choice  (2.7) gives
\bea
\lefteqn{\Int\frac{(bq^x,q^{1-x}\!/a;\:q)_\infty}{(-q^x,-q^{1-x};\:q)_\infty}\,
e^{\pi i x}(az)^x\,p(x)\,dx \label{c8} }\\
& & =\frac{(q,b/a;\:q)_\infty}{(z,b/az;\:q)_\infty}\Integ \frac{(azq^x,
q^{1-x}\!/az;\:q)_\infty}{(-q^x,-q^{1-x};\:q)_\infty}\,
(az)^x\,p(x)\,dx,\nonumber
\eea
where we  may assume that $p(x)$ is independent of $a,b$ and $z$. Let us
now replace $a,b,z$ by $-q^{-a},-q^b,-q^{a+c}$ in (\ref{c5}), respectively.
Then (\ref{c8}) can be written as
\bea
\lefteqn{\Int q^{cx}\,\frac{(-q^{b+x},-q^{a+1-x};\:q)_\infty}{(-q^x,
-q^{1-x};\:q)_\infty}\,p(x)\,dx \label{c9} }\\
& & =\frac{(q,q^{a+b};\:q)_\infty}{(q^{a+c},q^{b-c};\:q)_\infty}\,
\Integ q^{cx}\,\frac{(-q^{c+x},-q^{1-c-x};\:q)_\infty}{
(-q^x,-q^{1-x};\:q)_\infty}\,p(x)\,dx \nonumber
\eea
with \re$(a+c)>0$ and \re$(b-c)>0$. Setting $p(x)\equiv 1$ and denoting the
integral on the right side by $g(c)$ we may rewrite (\ref{c9}) in the form
\bea
\lefteqn{ g(c)=\frac{(q^{a+c},q^{b-c};\:q)_\infty}{(q,q^{a+b};\:q)_\infty}
\Int q^{cx}\,\frac{(-q^{b+x},-q^{a+1-x};\:q)_\infty}{
-q^x,-q^{1-x};\:q)_\infty}\,dx }\label{c10} \\ 
& & =\frac{1-q}{\log q^{-1}}\,\frac{\G{a+b}}{\G{a+c}\,\G{b-c}}\,
\Inte t^{c-1}\,\frac{(-tq^b,-q^{a+1}\!/t;\:q)_\infty}{
(-t,-q/t;\:q)_\infty}\,dt,\nonumber
\eea
assuming, without loss of generality, that $0<\re c<1$. \par

Since the left hand side is independent of $a,b$, we can set whatever
values of $a,b$ we wish, subject to the restriction mentioned above,
to compute $g(c)$. The simplest choice of $a,b$ is $a=0,\,b=1$. Then
\bea
\lefteqn{ g(c)=\frac{1-q}{\log q^{-1}}\,\frac{\G{1}}{\G{c}\G{1-c}}\,
\Inte t^{c-1}\,\frac{(-tq;\:q)_\infty}{(-t;\:q)_\infty}\,dt}\label{c11} \\
& & =\frac{1-q}{\log q^{-1}}\,\frac{\g{c}\,\g{1-c}}{\G{c}\,\G{1-c}}\,,\nonumber
\eea
since \dis
\[\Inte t^{c-1}\,\frac{(-tq;\:q)_\infty}{(-t;\:q)_\infty}\,dt=
\Inte\frac{t^{c-1}}{1+t}\,dt=\g{c}\g{1-c}. \]\dis
So, by (\ref{c9}) and (\ref{c11}) we have
\bea\lefteqn{
 (\log q^{-1})\,\Int q^{cx}\,\frac{(-q^{b+x},-q^{a+1-x};\:q)_\infty}{(-q^x,
-q^{1-x};\:q)_\infty}\,dx }\label{c12}   \\
& & =\frac{\g{c}\g{1-c}\G{a+c}\G{b-c}}{ \G{c}\G{1-c}\G{a+b} }\,. \nonumber
\eea
Substituting $q^x=t$ on the left and changing $a,b$ to $a-c$ and $b+c$,
respectively, we get (\ref{b7}).

It is also clear that setting $a=0$ in
(\ref{c12}) gives \r's formula (\ref{b3}). \par

We now explore a special choice for $p(x)$.
 We set $b=q^\a,\,a=q^{1-\beta},\,z=-q^{\beta-1/2},\:p(x)\equiv 1$ in
(\ref{c5}), getting
\be
\Int\frac{(q^{\a+x},q^{\beta-x};\:q)_\infty}{(-q^{1/2+x},-q^{1/2-x};
\:q)_\infty}\,dx=\frac{(q,q^{\a+\beta-1};\:q)_\infty}{
(-q^{\a-1/2},-q^{\beta-1/2};\:q)_\infty}, \label{c13}
\ee
\re$(\a+\beta-1)>0$, which is a $q$-analogue of yet another formula of
\r\ [22]:
\be
\Int\frac{dx}{\g{\a+x}\,\g{\beta-x}}=\frac{2^{\a+\beta-2}}{\g{\a+\beta-1}}
\,.\label{c14}
\ee
Using (\ref{c13}) one can show in a straightforward manner that
\bea
\lefteqn{\Int\frac{J_{\l+x}^{(1)}(a)\,J_{\mu-x}^{(2)}(b)}{a^{\l+x}\,b^{\mu-x}}
\,\frac{(-q^{\l+1/2},-q^{\mu+1/2};\:q)_\infty}{(-q^{1/2+x},-q^{1/2-x};
\:q)_\infty}\,dx}\label{c15} \\
& & =\frac{(q^{\l+\mu+1};\:q)_\infty}{(q;\:q)_\infty}\,\phy{2}{1}\left[
\ba{r}-q^{\l+1/2},\;-\frac{b^2}{a^2}q^{\mu+1/2}\\q^{\l+\mu+1}\ea;\:q,
-\frac{a^2}4
\right]\,,\nonumber
\eea
$|a/2|<1$, where
\be
J_{\nu}^{(1)}(x)=\sum_{m=0}^\infty\frac{(-1)^m(x/2)^{\nu+2m}}{(q;\:q)_m}\,\frac{
(q^{\nu+1+m};\:q)_\infty}{(q;\:q)_\infty}\,,
\ee\dis\dis\dis
\[
J_{\nu}^{(2)}(x)=\sum_{m=0}^\infty\frac{(-1)^m(x/2)^{\nu+2m}}{(q;\:q)_m}\,\frac{
(q^{\nu+1+m};\:q)_\infty}{(q;\:q)_\infty}\,q^{m(\nu+m)}\,,
\]
are the two $q$-Bessel functions of Jackson, see [15].  Formula (\ref{c15}) is
a $q$-analogue of the following formula 
\be
\Int\frac{J_{\l+x}(a)\,J_{\mu-x}(b)}{a^{\l+x}\,b^{\mu-x}}\,dx=\left(
\frac2{a^2+b^2}\right)^{\frac{\l+\mu}2}J_{\l+\mu}(\sqrt{2(a^2+b^2)}),
\ee
due to  \r\ in [21].

\section{ An integral Analogue of Bailey's \ph{6}{6} Sum.}
\setcounter{equation}{0} 
An application of (\ref{b11}) to (\ref{b22}) gives
\bea
J&=&\Integ(1-aq^{2x})(aq^{x+1}\!/b,aq^{x+1}\!/c,aq^{x+1}\!/d,aq^{x+1}\!/e;
\:q)_\infty \label{d1} \\
& & \mbox{}\cdot(q^{1-x}\!/b,q^{1-x}\!/c,q^{1-x}\!/d,q^{1-x}\!/e;
\:q)_\infty a^{2x}q^{2x^2-x}\omega(x) \nonumber \\
& & \mbox{}\hspace{-0.2in}\cdot\ph{6}{6}\left[\ba{c}q^{x+1}a^{1/2},-q^{x+1}a^{1/2},
\;\;bq^x,\;\;cq^x,\;\;dq^x,\;\;eq^x\\q^xa^{1/2},-q^xa^{1/2},aq^{x+1}\!/b,aq^{x+1}\!/c,
aq^{x+1}\!/d,aq^{x+1}\!/e\ea\,;\;q,\frac{qa^2}{bcde}\right]\,dx \nonumber \\
 &=&\frac{(aq/bc,aq/bd,aq/be,aq/cd,aq/ce,aq/de;\:q)_\infty}{
(qa^2/bcde;\:q)_\infty} \nonumber \\
& & \mbox{}\,\cdot\Integ(aq^{2x},q^{1-2x}\!/a;\:q)_\infty a^{2x}q^{2x^2-x}
\omega(x)\,dx, \nonumber
\eea
after we apply (\ref{b20}). Replacing $a$ by $\a^2$ and $q/b,q/c,q/d,q/e
$ by $a/\a,b/\a,c/\a,d/\a$, respectively, and taking
\be
\omega(s)=\frac{q^{s-2s^2}\a^{-4s}}{(\a^2q^{2s},q^{1-2s}\!/\a^2;\:q)_\infty}
p(s)\,,\quad p(s\pm1)=p(s),
\ee
we establish the relationship
\bea
\lefteqn{\Int\frac{(\a aq^x,aq^{-x}\!/\a,\a bq^x,bq^{-x}\!/\a,
\a cq^x,cq^{-x}\!/\a,\a dq^x,dq^{-x}\!/\a;\:q)_\infty}{
(\a^2q^{2x+1},q^{1-2x}\!/\a^2;\:q)_\infty}\,p(x)\,dx}\label{d3} \\
& & =\frac{(q,ab/q,ac/q,ad/q,bc/q,bd/q,cd/q;\:q)_\infty}{(abcd/q^3;\:q)_\infty}
\,\Integ p(x)\,dx, \nonumber
\eea
where $|abcd/q^3|<1$. In order to avoid singularities we also need to assume
that $\arg \a^2\neq 2k\pi,\;k=0,\pm1,\pm2,\ldots$. When $p(s)\equiv1$ and
$\a=i$ we obtain Askey's formula [4]:
\bea
\lefteqn{\Int \frac{h(i\sinh u;\:a,b,c,d)}{h(i\sinh u;q^{1/2},-q^{1/2},q
,-q)}\,du }\label{d4} \\
& & =(\log q^{-1})\frac{(ab/q,ac/q,ad/q,bc/q,bd/q,cd/q,q;\:q)_\infty}{
(abcd/q^3;\:q)_\infty}\,, \nonumber
\eea
where
\be h(x;\,a)=\prod_{n=0}^{\infty}(1-2axq^n+a^2q^{2n}),
\ee\dis\dis\dis
\[\mbox{}\quad h(x;\:a_1,a_2,\ldots,a_r)=\prod_{k=1}^{r}h(x;\,a_k). \]
%newly inserted

Askey's formula (\ref{d4}) was obtained from (\ref{d3}) by taking
$\a=i$ and specializing the unit periodic function $p(x)$ to be equal to 1.
But this is not the only case that we can evalute exactly. We will show
now that the integral $\Integ p(x)\,dx$ can be evaluated even in the
complicated case when $\a=i$ and
\be
p(x)=\frac{(-q^{2x+1},\;-q^{1-2x}\;q)_\infty\;\;(q^{-x}+q^x)}{
(fq^{-x},q^{x+1}/f,-fq^x,-q^{1-x}/f,gq^{-x},q^{x+1}/g,-gq^x,-q^{1-x}/g
;\;q)_\infty}\,, \label{n1} 
\ee
where \im $f$, \im$g$ and \im $(f/g)$ are not $0\!\!\!\pmod{2\pi}$ and
\im $(fg) \neq \pi \!\!\!\pmod{2\pi}$.  This will lead in a very natural way
to the Ismail-Masson q-beta integral (1.25). It can be easily verified that
$p(x\pm1)=p(x)$. Use of (\ref{d3}) then gives
\bea
\lefteqn{\eta(f,g)=\frac{(abcd/q^3;\;q)_\infty}{
(q,ab/q,ac/q,ad/q,bc/q,bd/q,cd/q;\;q)_\infty} }\label{n2} \\
& & \cdot \Int \frac{(iaq^x,-iaq^{-x},ibq^x,-ibq^{-x},icq^x,-icq^{-x},
idq^x,-idq^{-x};\;q)_\infty}{
(fq^{-x},q^{x+1}/f,-fq^x,-q^{1-x}/f,gq^{-x},q^{x+1}/g,-gq^x,-q^{1-x}/g
;\;q)_\infty} \nonumber \\
& & \mbox{\hspace{.2in}}\cdot (q^{-x}+q^x)\,dx\,, \nonumber
\eea
where
\be
\eta(f,g)=\Integ \frac{(-q^{2x+1},\;-q^{1-2x};\;q)_\infty\;\;(q^{-x}+q^x)\;dx}{
(fq^{-x},q^{x+1}/f,-fq^x,-q^{1-x}/f,gq^{-x},q^{x+1}/g,-gq^x,-q^{1-x}/g
;\;q)_\infty}\,. \label{n3}
\ee
Observe that $\eta(f,g)$ is independent of $a,b,c,d$, so the expression on
the right hand side of (\ref{n2}) must have the same property. For the
purpose of (\ref{n2}) we then set
\be
ab=q^2=cd,\qquad a=-iq/f,\quad c=iq \label{n4}
\ee
to get
\bea
\lefteqn{\eta(f,g)=\frac1{(q,q,g/f,fq^2/g,-fg,-q^2/fg;\;q)_\infty}
}\label{n5}\\
& & \cdot\Int \frac{(q^{-x}+q^x)\,dx}{[1+f(q^x-q^{-x})-f^2][1-q/g(q^x-
q^{-x})-q^2/g^2]}. \nonumber
\eea
Substituting $q^{-x}-q^x=u$, the integral on the right side of (\ref{n5})
becomes
\bea\lefteqn{
\frac1{\log q^{-1}}\Int\frac{du}{(1-f^2-fu)(1-q^2/g^2+qu/g) } }\label{n6} \\
& & =\frac1{f\,\log q^{-1}}\;\frac{2\pi i}{(1-qf/g)(1+q/fg)}. \nonumber
\eea
Formula (3.11) can be proved by either using  a partial fraction decomposition or by a simple  contour integration. Thus
\be
\eta(f,g)=\frac{2\pi i}{f\,\log q^{-1}\,(q,q,g/f,qf/g,-fg,-q/fg;\;q)_\infty}
\label{n7}
\ee
Combining (\ref{d3}), (\ref{n2}), (\ref{n3}) and (\ref{n7}) we find that
\bea\lefteqn{
\Int\frac{(iaq^x,-iaq^{-x},ibq^x,-ibq^{-x},icq^x,-icq^{-x},
idq^x,-idq^{-x};\;q)_\infty \,(q^{-x}+q^x)}{
(fq^{-x},q^{x+1}/f,-fq^x,-q^{1-x}/f,gq^{-x},q^{x+1}/g,-gq^x,-q^{1-x}/g
;\;q)_\infty}\,dx }\label{n8} \\
& & =\frac{2 \pi i\;(ab/q,qc/q,ad/q,bc/q,bd/q,cd/q;\;q)_\infty}{
f\,\log q^{-1}\,(q,g/f,qf/g,-fg,-q/fg,abcd/q^3;\;q)_\infty}.\nonumber
\eea
This is the same as (1.25) which was established in [16] by an
entirely different method.

We would like to mention that formula (\ref{d3}) is valid even when $\a$ is
real, provided the integral on the left is interpreted as a principal-value
integral. For a detailed discussion of this point, see [18]. \par

We should like to point out that if we denote
\be
f(t)=\frac{(\a at,a/\a t,\a bt,b/\a t,\a ct,c/\a t,\a dt,d/\a t;\:q)_\infty}{
(q\a^2t^2,q/\a^2t^2;\:q)_\infty\,t\,(\log q^{-1})}\,,
\ee
then (\ref{b20}) and (\ref{d3}) state that 
\bea
\lefteqn{\Inte f(t)\,dt=\Inte \left(\frac{\log q^{-1}}{1-q}\right)\,f(t)\,
d_qt }\label{d7} \\
& & =\frac{(q,ab/q,ac/q,ad/q,bc/q,bd/q,cd/q;\:q)_\infty}{(abcd/q^3;\:q)_\infty}
\,,\nonumber
\eea
where
\be
\Inte g(t)\,d_qt=(1-q)\sum_{n=-\infty}^{\infty}g(q^n)q^n
\ee
is a $q$-integral defined by Jackson, see [11]. This is an example where an
absolutely continuous measure and a purely discrete measure on the real line
have the same total weight, indicating an indeterminate moment problem.
Ismail and Masson [16] have recently found two systems of rational functions
which are biorthogonal with respect to the weight function given in the integral
(1.23).  In earlier unpublished work,  Rahman proved the biorthogonality of the same rational functions with respect to the weight function in the 
Askey integral (\ref{d4}). Further properties of these rational functions, 
their Rodrigues formulas, the $q$-difference equations etc. will be discussed 
in a subsequent paper. For a different system of biorthogonal rational 
functions on $[-1,1]$, see [19].

\section{ Case of an anti-unit-periodic function.}
\setcounter{equation}{0} 
We shall now derive a $q$-analogue of \r's formula [22]:
\bea
\Int \frac{p(x)\,dx}{\g{\a+x}\g{\beta-x}\g{\gamma+x}\g{\delta-x}}\label{e1} \\
= \frac{\Integ p(t)\cos [\pi(2t+\a-\beta)/2]\,dt}{\g{\a+\beta/2}\g{\gamma+
\delta/2} \g{\a+\delta-1}},\nonumber
\eea
where $p(x\pm1)=-p(x),\;\a+\delta=\beta+\gamma$ and $\re (\a+\beta+\gamma+
\delta)>2$. Let us first write down the \ph{2}{2} summation formula
[11, (5.3.4)] in the following form:
\bea
\lefteqn{\sum_{n=-\infty}^\infty(aq^{n+1}\!/b,aq^{n+1}\!/c,q^{1-n}\!/b,
q^{1-n}\!/c;\:q)_\infty a^ne^{\pi in}q^{n^2} } \label{e2} \\
& & =\frac{(aq/bc;\:q)_\infty}{(-aq/bc;\:q)_\infty}\,(q^2,aq,q/a,aq^2/b^2,
aq^2/c^2;\:q^2)_\infty,\quad |aq/bc|<1. \nonumber
\eea
Let $\omega(x)$ be a bounded continuous function on {\bf R} such that
$\omega(x\pm1)=-\omega(x)$. Then it can be shown that
\bea
\lefteqn{\Int (aq^{x+1}\!/b,aq^{x+1}\!/c,q^{1-x}\!/b,q^{1-x}\!/c;\:q)_\infty
a^xq^{x^2}\omega(x)\,dx} \label{e3} \\
& & =\Integ(aq^{x+1}\!/b,aq^{x+1}\!/c,q^{1-x}\!/b,q^{1-x}\!/c;\:q)_\infty
a^xq^{x^2}\omega(x) \nonumber \\
& & \mbox{}\quad\cdot\ph{2}{2}\left[\ba{cc}bq^x,&cq^x\\aq^{x+1}\!/b,
& aq^{x+1}\!/c\ea\,;\;q,-aq/bc\right]\,dx\nonumber \\
& & =\frac{(aq/bc;\,q)_\infty}{(-aq/bc;\,q)_\infty}\,(q^2,aq^2/b^2,aq^2/c^2;
\:q^2)_\infty 
\nonumber \\& & \mbox{}\quad
\Integ(aq^{2x+1},q^{1-2x}\!/a;\:q^2)_\infty a^2q^{x^2}
\,\omega(x)\,dx\,. \nonumber
\eea
Replacing $a,b,c$ by $q^{\a-\beta},q^{1-\beta},q^{1-\delta}$, respectively,
 and assuming the $\a+\delta=\beta+\gamma$, this can be written as
\bea
\lefteqn{\Int(q^{\a+x},q^{\beta-x},q^{\gamma+x},q^{\delta-x};\:q)_\infty 
q^{(\a-\beta)x+x^2}\omega(x)\,dx}\label{e4} \\
& & =\frac{(q^{\a+\delta-1};\:q)_\infty}{(-q^{\a+\delta-1};\:q)_\infty}
(q^{\a+\beta},q^{\gamma+\delta},q^2;\:q^2)_\infty \nonumber \\
& & \mbox{}\quad\cdot\Integ(q^{\a-\beta+1+2x},q^{\beta-\a+1-2x};\:q^2)_\infty
q^{(\a-\beta)x+x^2}\omega(x)\,dx\,. \nonumber 
\eea
In terms of the $q$-gamma function this can be written in the form
\bea
\lefteqn{\Int \frac{q^{(\a-\beta)x+x^2}\omega(x)\,dx\,}{\G{\a+x}\G{\beta-x}
\G{\gamma+x}\G{\delta-x}}}\label{e5} \\
& & =\frac{(-q;\,q)_{\a+\delta-2}}{(1+q)^{\a+\delta-2}}\,\frac{
\Gamma_{q^2}^2(\frac12)}{\Ga{\frac{\a+\beta}2}\Ga{\frac{\gamma+\delta}2}
\G{\a+\delta-1 }}\nonumber \\
& & \mbox{}\quad\cdot\Integ\frac{q^{(\a-\beta)x+x^2}\omega(x)\,dx\,}{\Ga{
\frac{\a-\beta+1}2+x}\,\Ga{\frac{\beta-\a+1}2-x} }\,, \nonumber
\eea
where\dis
\be
(a;\,q)_\l=\frac{(a;\,q)_\infty}{(aq^\l;\,q)_\infty}\,,
\ee
$\a+\delta=\beta+\gamma$ and \re$(\a+\beta+\gamma+\delta)>2$.
It is easy to see that in the limit $q\to 1^-$ formula (\ref{e5}) approaches
(\ref{e1}).

\section{ An integral analogue of the \ph{8}{8} sum.}
\setcounter{equation}{0} 
It is clear from (\ref{b11}) that we can associate a general \ph{r}{r} series
with a doubly infinite integral of the type considered in \S3. However,
there are no known summation formulas for a very-well-poised \ph{r}{r} for
$r>6$, only transformation formulas. One could use these transformation
formulas, see [11, chapter 5], to express such integrals in terms of a
string of basic hypergeometric series, but the exercise does not seem to
have a purpose, except for the case $r=8$. The results that we shall obtain
in this section will, hopefully, convince the reader that the case of  \ph{8}{8}
series has some interesting features. \par

Let us consider the integral
\bea
\lefteqn{ K:=\Int(1-aq^{2x})(aq^{x+1}\!/b,aq^{x+1}\!/c,aq^{x+1}\!/d,
aq^{x+1}\!/e,aq^{x+1}\!/f;\:q)_\infty \label{f1} }\\
& & \mbox{}\quad\cdot\frac{(q^{1-x}\!/b,q^{1-x}\!/c,q^{1-x}\!/d,q^{1-x}\!/e,
q^{1-x}\!/f;\:q)_\infty}{(aq^{x+1}\!/g,q^{1-x}\!/g;\:q)_\infty}\,
q^{2x^2-x}a^{2x}\omega(x)\,dx\,,\nonumber
\eea
where $\omega(x)$ has the same properties as mentioned in \S3. Notice that
we have taken a pair of infinite products in the denominator also, which 
makes the structure of $K$ slightly different from that of  the integral $J$ defined in
(\ref{b22}) and (\ref{d1}). By exploiting the unit-periodic property of
$\omega(x)$ we could avoid this difference but we believe the form of the 
integrand in (\ref{f1}) is more instructive. It is obvious that when $g$
equals any one of the parameters $b,c,d,e,f$ in the numerator then $K$
will reduce to $J$. We will assume that
\be
\left|\frac{qga^2}{bcdef}\right|<1\,,
\ee
which ensures the convergence of the integral. An application of (\ref{b11}) then gives
\be
K= \label{f3} 
\ee\dis\dis\dis
\bea
\lefteqn{
\Integ( 1-aq^{2x})(aq^{x+1}\!/b,aq^{x+1}\!/c,aq^{x+1}\!/d,
aq^{x+1}\!/e,aq^{x+1}\!/f;q)_\infty }\nonumber \\
& & \mbox{}\cdot\frac{(q^{1-x}\!/b,q^{1-x}\!/c,q^{1-x}\!/d,q^{1-x}\!/e,
q^{1-x}\!/f;\:q)_\infty}{(aq^{x+1}\!/g,q^{1-x}\!/g;\:q)_\infty}\,
q^{2x^2-x}\omega(x)\nonumber \\
& & \mbox{}\hspace{-0.6in}\cdot\ph{8}{8}\left[\ba{c}q^{x+1}a^{\frac12},
-q^{x+1}a^{\frac12},
\;\;bq^x,\;\;cq^x,\;\;dq^x,\;\;eq^x,\;\;fq^x,aq^{x+1}\!/g \\ q^xa^{\frac12},
-q^xa^{\frac12},aq^{x+1}\!/b,aq^{x+1}\!/c,aq^{x+1}\!/d,aq^{x+1}\!/e,
aq^{x+1}\!/f,gq^x  \ea ;\:q,\frac{qga^2}{
bcdef}\right]\,dx\,. \nonumber
\eea
Using the transformation formula [11,(5.6.2)] we get, for the \ph{8}{8} 
series above
\bea
\lefteqn{\ph{8}{8}[\;]=\frac{(q,aq/bf, aq/cf,aq/df, aq/ef,qf/b,qf/c,qf/d,qf/e;
\:q)_\infty}{
(qf/g,aq/fg,qf^2/a;\:q)_\infty(aq^{x+1}\!/b,aq^{x+1}\!/c,aq^{x+1}\!/d,
aq^{x+1}\!/e;\:q)_\infty} }\label{f4} \\
& & \mbox{} \cdot\frac{(aq^{x+1}\!/g,q^{1-x}\!/g,aq^{2x+1},q^{-2x}\!/a;\:q)_\infty}{
(aq^{x+1}\!/f, q^{1-x}\!/b,q^{1-x}\!/c,q^{1-x}\!/d,q^{1-x}\!/e,
q^{1-x}\!/f;\:q)_\infty}\nonumber \\
& & \mbox{}\quad\cdot \mbox{}_8W_7(f^2/a;\:bf/a,cf/a,df/a,ef/a,gf/a;\:q,qga^2/bcdf)\nonumber \\
& &\mbox{}\hspace{-0.1in}  +\frac{(q, aq^2/bg,aq^2/cg,aq^2/dg,aq^2/eg,g/b,g/c,g/d,g/e;\;q)_\infty}{
(qf/g,fg/aq,aq^3/g^2;\,q)_\infty( aq^{x+1}\!/b,aq^{x+1}\!/c,aq^{x+1}\!/d,
aq^{x+1}\!/e;\:q)_\infty} \nonumber  \\
& & \mbox{} \cdot\frac{(fq^x,fq^{-x}\!/a,aq^{2x+1},q^{-2x}\!/a;\:q)_\infty}
{(gq^x,gq^{-x}\!/a,q^{1-x}\!/b,q^{1-x}\!/c,q^{1-x}\!/d,q^{1-x}\!/e
;\:q)_\infty}\nonumber \\
& & \mbox{}\cdot \mbox{}_8W_7(aq^2/g^2;\:bq/g,cq/g,dq/g,eq/g,fq/g
;\:q, qga^2/bcdf),\nonumber 
\eea
where
\bea
\lefteqn{
\mbox{}_8W_7(a;\: b,c,d,e,f;\:q,z)}\label{f5} \\
& & \mbox{}\hspace{-0.2in}:=\phy{8}{7}\left[\ba{c}\!a,qa^{\frac12},-qa^{\frac12},\quad b,\quad\;c,
\quad\;d, \quad\;e,\quad\;f\\
\mbox{}\quad\;a^{\frac12},\;-a^{\frac12},aq/b,aq/c,aq/d,aq/e,aq/f\ea\,;\;q,
z\right]. \nonumber
\eea
Choosing
\be
\omega(s)=\frac{q^{s-2s^2}a^{-2s}p(s)}{(aq^{2s},q^{-2s}\!/a;\:q)_\infty}\,,
\qquad p(s\pm1)=p(s),
\ee
in (5.3) we obtain the following identity
\bea
\lefteqn{\Int\frac{( aq^{x+1}\!/b,aq^{x+1}\!/c,aq^{x+1}\!/d,
aq^{x+1}\!/e,aq^{x+1}\!/f;\:q)_\infty}{(aq^{2x+1},q^{1-2x}\!/a;\:q)_\infty}
 \label{f7} }\\
& & \mbox{}\cdot\frac{(q^{1-x}\!/b,q^{1-x}\!/c,q^{1-x}\!/d,q^{1-x}\!/e,
q^{1-x}\!/f;\:q)_\infty}{(aq^{x+1}\!/g,q^{1-x}\!/g;\:q)_\infty}\,
p(x)\,dx\nonumber \\
& & = \frac{(q,aq/bf, aq/cf,aq/df, aq/ef,qf/b,qf/c,qf/d,qf/e;\:q)_\infty}{
(aq/fg,qf/g, qf^2/a;\:q)_\infty} \nonumber\\
& & \mbox{}\hspace{-0.1in}\cdot \mbox{}_8W_7(f^2/a;\:bf/a,cf/a,df/a,ef/a,gf/a;\:q,
qga^2/bcdf)\cdot\Integ p(x)\,dx\nonumber \\
& &\mbox{}+\frac{(q,g/b,g/c,g/d,g/e , aq^2/bg,aq^2/cg,aq^2/dg,
aq^2/eg;\:q)_\infty}{ (fg/aq,aq^3/g^2,qf/g ;\,q)_\infty} \nonumber  \\
& & \mbox{}\cdot \mbox{}_8W_7(aq^2/g^2;\:bq/g,cq/g,dq/g,eq/g,fq/g
;\:q, qga^2/bcdf)\nonumber  \\
& & \mbox{}\cdot\Integ\frac{(fq^x,q^{1-x}\!/f,aq^{x+1}\!/f,
fq^{-x}\!/a;\:q)_\infty}{(gq^x,q^{1-x}\!/g,aq^{x+1}\!/g,gq^{-x}\!/a;
\:q)_\infty}\,p(x)\,dx, \nonumber
\eea
provided $\arg a$ and $\arg g$ are neither 0 nor multiples of $2\pi$. \par

Replacing $a$ by $\a^2$ and $q/b,q/c,q/d,q/e,q/f,q/g$ by $a/\a,
b/\a,c/\a$,$d/\a$,$e/\a$ and $f/\a$, respectively, we obtain
\bea
\lefteqn{\Int \frac{(a\a q^x,a q^{-x}\!/\a,b\a q^x,bq^{-x}\!/\a,
c\a q^x,cq^{-x}\!/\a,
d\a q^x,dq^{-x}\!/\a;\:q)_\infty}{(\a^2q^{2x+1},q^{1-2x}\!/\a^2;\:q)_\infty
 }}\label{f8} \\
& & \mbox\quad\cdot\frac{(e\a q^x,eq^{-x}\!/\a;\:q)_\infty}{(f\a q^x,
fq^{-x}\!/\a;\:q)_\infty}\,p(x)\,dx \nonumber \\
& & =\frac{(q,aq/e,bq/e,cq/e,dq/e,ae/q,be/q,ce/q,de/q;\:q)_\infty}{
(ef/q,qf/e,q^3/e^2;\:q)_\infty} \nonumber \\
& & \mbox{}\hspace{-0.0in}\cdot \mbox{}_8W_7(q^2/e^2;\:q^2/ae,q^2/be,q^2/ce,q^2/de,
qf/e ;\:q, abcde/fq^3)\Integ p(x)\,dx\nonumber  \\
& & +\frac{(q,af,bf,cf,df,q/f,b/f,c/f,d/f;\:q)_\infty}{(q/ef,qf/e,qf^2;\:
q)_\infty} \nonumber \\
& & \cdot \mbox{}_8W_7(f^2;\:qf/a,qf/b,qf/c,qf/d,qf/e;
\:q, abcde/fq^3)\nonumber  \\
& & \cdot\Integ\frac{(\a q^{x+1}\!/e,eq^{-x}\!/\a,\a eq^{x+1},
q^{-x}\!/\a e;\:q)_\infty} {(\a q^{x+1}\!/f,fq^{-x}\!/\a,\a fq^{x+1},
q^{-x}\!/\a f;\:q)_\infty} \,p(x)\,dx\,, \nonumber
\eea
where the parameters $\a$ and $f$ are such that no zeros occur in any
denominator. Recall the transformation formula [11, (2.11.1)],
\be
\mbox{}_8W_7(q^2/e^2;\:qf/e,q^2/ae,q^2/be,q^2/ce,q^2/de
;\:q, abcde/fq^3) \label{f9} 
\ee\dis\dis
\bea
\lefteqn{ 
 =\frac{(q^3/e^2,cd/q,ac/q,ad/q,bq/d,bq/a,ef/q,fq^3/ade;\:q)_\infty}{
(aq/e,cq/e,dq/e,be/q,acde/q^3,bq^3/ade,fq/d,fq/a;\:q)_\infty} \nonumber }\\
& & \cdot \mbox{}_8W_7(bq^2/ade;\:b/f,bc/q,q^2/ad,q^2/de,q^2/ae;
\:q,fq/c) \nonumber  \\
& &\mbox{}\hspace{-0.2in}-\frac{(q^3/e^2,ef/q,bf,cf,df,af,q^2/ae,q^2/ce,q^2/de,b/f,acd/fq^2,
fq^3/acd;\:q)_\infty}{(q/ef,aq/e,bq/e,cq/e,dq/e,be/q,fq/a,fq/c,fq/d,acde/q^3,
q^4/acde,qf^2;\:q)_\infty}\nonumber \\
& & \cdot \mbox{}_8W_7(f^2;\:qf/a,qf/b,qf/c,qf/d,qf/e
;\:q, abcde/fq^3)\,. \nonumber 
\eea
Note that the last $\mbox{}_8W_7$ series on the right is the same as that
on the right side of (\ref{f8}). This enables us to rewrite (\ref{f8}) in
the form
\be \mbox{}\label{f10}
\ee\dis\dis\dis
\bea
\lefteqn{\mbox{}\hspace{-0.4in} \Int \frac{(a\a q^x\!,a q^{-x}\!/\a,b\a q^x\!,
bq^{-x}\!/\a,c\a q^x\!,cq^{-x}\!/\a,
d\a q^x\!,dq^{-x}\!/\a,e\a q^x\!,eq^{-x}\!/\a;q)_\infty}{(\a^2q^{2x+1},
q^{1-2x}\!/\a^2,f\a q^x,fq^{-x}\!/\a;q)_\infty }p(x)dx}\nonumber\\
& & \mbox{}\hspace{-0.4in} =\frac{(q,ac/q,ad/q,cd/q,ae/q,ce/q,de/q,bq/a,bq/d,bq/e,fq^3/ade;
\:q)_\infty}{(qf/a,qf/d,qf/e,bq^3/ade,acde/q^3;\:q)_\infty} \nonumber \\
& & \cdot \mbox{}_8W_7(bq^2/ade;\:b/f,bc/q,q^2/ad,q^2/de,q^2/ae;
\:q,fq/c)  \Integ p(x)\,dx\nonumber  \\
& & + \; \mbox{}_8W_7(f^2;\:qf/a,qf/b,qf/c,qf/d,qf/e
;\:q, abcde/fq^3) \nonumber  \\
& &\mbox{}\;\cdot\left\{
\frac{(q,a/f,b/f,c/f,d/f,af,bf,cf,df;\:q)_\infty}{(qf/e,q/ef,qf^2;\:q)_\infty}
\,\l(e,f)\nonumber\right. \\
& &\mbox{\hspace{.2in}}- \frac{(q,af,bf,cf,df,b/f,ae/q,ce/q,de/q,q^2/ae,q^2/ce,q^2/de;
\:q)_\infty}{(qf/e,q/ef,qf^2,fq/a,fq/c,fq/d;\:q)_\infty}\nonumber \\
& & \mbox{} \left. \quad\cdot\frac{(acd/fq^2,fq^3/acd;\:q)_\infty}{(acde/q^3,q^4/acde;
\:q)_\infty}\,\Integ p(x)\,dx\right\}\,,\nonumber
\eea
where
\be
\l(e,f)=\Integ\frac{(\a q^{x+1}\!/e, eq^{-x}\!/\a,\a eq^{x+1},
q^{-x}\!/\a e;\:q)_\infty} {(\a q^{x+1}\!/f,fq^{-x}\!/\a,\a fq^{x+1},
q^{-x}\!/\a f;\:q)_\infty} \,p(x)\,dx\,.
\ee
The special case $f=abcde/q^4$ is of particular interest because it gives
an overall balance of the parameters inside the integral on the left side
while reducing the first $\mbox{}_8W_7$ series on the right to a very-%
well-poised \phy{6}{5}, which is summable by use of [11, (2.7.1)]. This
leads to the formula
\pagebreak[1]
\bea\lefteqn{\mbox{}\hspace{-0.1in} \Int \frac{(a\a q^x,a q^{-x}\!/\a,b\a q^x,bq^{-x}\!/\a,c\a q^x,cq^{-x}\!/\a,
d\a q^x,dq^{-x}\!/\a,e\a q^x,eq^{-x}\!/\a,;\:q)_\infty}{(\a^2q^{2x+1},
q^{1-2x}\!/\a^2,\alpha abcde q^{x-4},abcdeq^{-x-4}\!/\a;\:q)_\infty }\,p(x)\,dx }\label{f12}
\\
& &\mbox{}\hspace{-0.3in} =\frac{(q,ab/q,ac/q,ad/q,ae/q,bc/q,bd/q,be/q,cd/q,ce/q,de/q; q)_\infty}
{(qf/a,qf/b,qf/c,qf/d,qf/e;\:q)_\infty}\Integ p(x)\,dx\nonumber \\
& &\mbox{}\hspace{-0.1in} +\frac{(q,af,bf,cf,df,a/f,b/f,c/f,d/f;\:q)_\infty}{
(qf^2,q/ef,qf/e;\:q)_\infty}\nonumber \\
& &\mbox{}\hspace{-0.2in}\cdot\left\{\l(e,f)-\frac{(ae/q,be/q,ce/q,de/q,q^2/ae,q^2/be,q^2/ce,q^2/de;
\:q)_\infty}{ (a/f,b/f,c/f,d/f,qf/a,qf/c,qf/c,qf/d;\:q)_\infty}
\Integ p(x)\,dx\right\} \nonumber \\
& & \cdot\, \mbox{}_8W_7(f^2;\:qf/a,qf/b,qf/c,qf/d,qf/e
;\:q, q)\,, \nonumber 
\eea
with 
\be
f=abcde/q^4\,.
\ee
This formula is essentially the same as the one found in [18].

\end{document}